\documentclass{amsart}
\usepackage{amsmath}
\usepackage{amsfonts}

\newtheorem{theorem}{Theorem}
\theoremstyle{plain}

\newtheorem{corollary}{Corollary}

\numberwithin{equation}{section}
\input{tcilatex}

\begin{document}
\title[\textbf{Kim's }$q$\textbf{-Euler numbers and polynomials }]{EXPLICIT
FORMULAS INVOLVING $q$-EULER NUMBERS AND POLYNOMIALS}
\author[\textbf{S. Arac\i }]{\textbf{Serkan Arac\i }}
\address{University of Gaziantep, Faculty of Science and Arts, Department of
Mathematics, 27310 Gaziantep, TURKEY}
\email{mtsrkn@hotmail.com}
\author[\textbf{M. Acikgoz}]{\textbf{Mehmet Acikgoz}}
\address{University of Gaziantep, Faculty of Science and Arts, Department of
Mathematics, 27310 Gaziantep, TURKEY}
\email{acikgoz@gantep.edu.tr}
\author[\textbf{J J. Seo}]{\textbf{Jong Jin Seo}}
\address{Department of Applied Mathematics, Pukyong National University,
Busan{} 608-737, Republic of Korea}
\email{seo2011@pknu.ac.kr (\textbf{Corresponding author})}
\subjclass[2000]{\textbf{\ Primary 05A10, 11B65; Secondary 11B68, 11B73.}}
\keywords{ \textbf{Euler numbers and polynomials, }$q$\textbf{-Euler numbers
and polynomials with weight }$0$\textbf{, }$q$\textbf{-Bernoulli numbers
with weight }$0$\textbf{, }$p$\textbf{-adic }$q$\textbf{-integral}}

\begin{abstract}
In this paper, we deal with q-Euler numbers and q-Bernoulli numbers. We
derive some interesting relations for q-Euler numbers and polynomials by
using their generating function and derivative operator. Also, we show
between the q-Euler numbers and q-Bernoulli numbers via the p-adic
q-integral in the p-adic integer ring.
\end{abstract}

\maketitle

\section{\textbf{PRELIMINARIES}}

\bigskip Imagine that $p$ be a fixed odd prime number. Throughout this paper
we use the following notations, by $%
\mathbb{Z}
_{p}$ denotes the ring of $p$-adic rational integers, $%
\mathbb{Q}
$ denotes the field of rational numbers, $%
\mathbb{Q}
_{p}$ denotes the field of $p$-adic rational numbers, and $%
\mathbb{C}
_{p}$ denotes the completion of algebraic closure of $%
\mathbb{Q}
_{p}$. Let $%
\mathbb{N}
$ be the set of natural numbers and $%
\mathbb{N}
^{\ast }=%
\mathbb{N}
\cup \left\{ 0\right\} $.

The $p$-adic absolute value is defined by 
\begin{equation*}
\left\vert p\right\vert _{p}=\frac{1}{p}.
\end{equation*}

In this paper we assume $\left\vert q-1\right\vert _{p}<1$ as an
indeterminate.

$\left[ x\right] _{q}$ is a $q$-extension of $x$ which is defined by 
\begin{equation*}
\left[ x\right] _{q}=\frac{1-q^{x}}{1-q},
\end{equation*}%
we note that $\lim_{q\rightarrow 1}\left[ x\right] _{q}=x$\ (see[1-12]).

We say that $f$ is a uniformly differentiable function at a point $a\in 
\mathbb{Z}
_{p}$, if the difference quotient 
\begin{equation*}
F_{f}\left( x,y\right) =\frac{f\left( x\right) -f\left( y\right) }{x-y}
\end{equation*}

has a limit $f%
{\acute{}}%
\left( a\right) $ as $\left( x,y\right) \rightarrow \left( a,a\right) $ and
denote this by $f\in UD\left( 
\mathbb{Z}
_{p}\right) $.

\bigskip Let $UD\left( 
\mathbb{Z}
_{p}\right) $ be the set of uniformly differentiable function on $%
\mathbb{Z}
_{p}$. For $f\in UD\left( 
\mathbb{Z}
_{p}\right) $, let us start with the expressions%
\begin{equation*}
\frac{1}{\left[ p^{N}\right] _{q}}\sum_{0\leq \xi <p^{N}}f\left( \xi \right)
q^{\xi }=\sum_{0\leq \xi <p^{N}}f\left( \xi \right) \mu _{q}\left( \xi +p^{N}%
\mathbb{Z}
_{p}\right) ,
\end{equation*}

represents $p$-adic $q$-analogue of Riemann sums for $f$. \ The integral of $%
f$ on $%
\mathbb{Z}
_{p}$ will be defined as the limit $\left( N\rightarrow \infty \right) $ of
these sums, when it exists. The $p$-adic $q$-integral of function $f\in
UD\left( 
\mathbb{Z}
_{p}\right) $ is defined by T. Kim 
\begin{equation}
I_{q}\left( f\right) =\int_{%
\mathbb{Z}
_{p}}f\left( \xi \right) d\mu _{q}\left( \xi \right) =\lim_{N\rightarrow
\infty }\frac{1}{\left[ p^{N}\right] _{q}}\sum_{\xi =0}^{p^{N}-1}f\left( \xi
\right) q^{\xi }\text{ }  \label{equation 1}
\end{equation}

The bosonic integral is considered as a bosonic limit $q\rightarrow 1,$ $%
I_{1}\left( f\right) =\lim_{q\rightarrow 1}I_{q}\left( f\right) $.
Similarly, the fermionic $p$-adic integral on $%
\mathbb{Z}
_{p}$ is introduced by T. Kim as follows:%
\begin{equation}
I_{-q}\left( f\right) =\lim_{q\rightarrow -q}I_{q}\left( f\right) =\int_{%
\mathbb{Z}
_{p}}f\left( \xi \right) d\mu _{-q}\left( \xi \right)  \label{equation 2}
\end{equation}

(for more details, see [9-12]).

In \cite{Kim 4}, the $q$-Euler polynomials with wegiht $0$ are introduced as%
\begin{equation}
\widetilde{E}_{n,q}\left( x\right) =\int_{%
\mathbb{Z}
_{p}}\left( x+y\right) ^{n}d\mu _{-q}\left( y\right)  \label{equation 108}
\end{equation}

From (\ref{equation 108}), we have%
\begin{equation*}
\widetilde{E}_{n,q}\left( x\right) =\sum_{l=0}^{n}\binom{n}{l}x^{l}%
\widetilde{E}_{n-l,q}
\end{equation*}

where $\widetilde{E}_{n,q}(0)=\widetilde{E}_{n,q}$ are called $q$-Euler
numbers with weight $0$. Then, $q$-Euler numbers are defined as%
\begin{equation*}
q\left( \widetilde{E}_{q}+1\right) ^{n}+\widetilde{E}_{n,q}=\left\{ \QATOP{%
\left[ 2\right] _{q},\text{ if }n=0}{0,\text{ \ \ \ \ if\ }n\neq 0,}\right.
\end{equation*}

with the usual convention about replacing $\left( \widetilde{E}_{q}\right)
^{n}$ by $\widetilde{E}_{n,q}$ is used.

Similarly, the $q$-Bernoulli polynomials and numbers with weight $0$ are
defined, respectively%
\begin{eqnarray*}
\widetilde{B}_{n,q}\left( x\right) &=&\lim_{n\rightarrow \infty }\frac{1}{%
\left[ p^{n}\right] _{q}}\sum_{y=0}^{p^{n}-1}\left( x+y\right) ^{n}q^{y} \\
&=&\int_{%
\mathbb{Z}
_{p}}\left( x+y\right) ^{n}d\mu _{q}\left( y\right)
\end{eqnarray*}

and 
\begin{equation*}
\widetilde{B}_{n,q}=\int_{%
\mathbb{Z}
_{p}}y^{n}d\mu _{q}\left( y\right)
\end{equation*}

(for more informations, see \cite{Kim 2}).

We, by using Kim's et al. method in \cite{Kim}, will investigate some
interesting identities on the $q$-Euler numbers and polynomials from their
generating function and derivative operator. Consequently, we derive some
properties on $q$-Euler numbers and polynomials and $q$-Bernoulli numbers
and polynomials by using $q$-Volkenborn integral and fermionic $p$-adic $q$%
-integral on $%
\mathbb{Z}
_{p}$.

\section{\textbf{ON KIM'S }$q$\textbf{-EULER NUMBERS AND POLYNOMIALS}}

Let us consider Kim's $q$-Euler polynomials as follows: 
\begin{equation}
F_{x}^{q}=F_{x}^{q}\left( t\right) =\frac{\left[ 2\right] _{q}}{qe^{t}+1}%
e^{xt}=\sum_{n=0}^{\infty }\widetilde{E}_{n,q}\left( x\right) \frac{t^{n}}{n!%
}.  \label{equation 100}
\end{equation}

Here $x$ is a fixed parameter. Thus, by expression of (\ref{equation 100}),
we can readily see the following 
\begin{equation}
qe^{t}F_{x}^{q}+F_{x}^{q}=\left[ 2\right] _{q}e^{xt}.  \label{equation 101}
\end{equation}

Last from equality, taking derivative operator $D$ as $D=\frac{d}{dt}$ on
the both sides of (\ref{equation 101}). Then, we easily see that%
\begin{equation}
qe^{t}\left( D+I\right) ^{k}F_{x}^{q}+D^{k}F_{x}^{q}=\left[ 2\right]
_{q}x^{k}e^{xt}  \label{equation 106}
\end{equation}

where $k\in 
\mathbb{N}
^{\ast }$ and $I$ is identity operator. By multiplying $e^{-t}$ on both
sides of (\ref{equation 106}), we get%
\begin{equation}
q\left( D+I\right) ^{k}F_{x}^{q}+e^{-t}D^{k}F_{x}^{q}=\left[ 2\right]
_{q}x^{k}e^{\left( x-1\right) t}  \label{equation 3}
\end{equation}

Let us take derivative operator $D^{m}\left( m\in 
\mathbb{N}
\right) $ on both sides of (\ref{equation 3}). Then we get 
\begin{equation}
qe^{t}D^{m}\left( D+I\right) ^{k}F_{x}^{q}+D^{k}\left( D-I\right)
^{m}F_{x}^{q}=\left[ 2\right] _{q}x^{k}\left( x-1\right) ^{m}e^{xt}
\label{equation 107}
\end{equation}

Let $G\left[ 0\right] $ (not $G\left( 0\right) $) be the constant term in a
Laurent series of $G\left( t\right) $. Then, from (\ref{equation 107}), we
get%
\begin{equation}
\sum_{j=0}^{k}\binom{k}{j}\left( qe^{t}D^{k+m-j}F_{x}^{q}\left( t\right)
\right) \left[ 0\right] +\sum_{j=0}^{m}\binom{m}{j}\left( -1\right)
^{j}\left( D^{k+m-j}F_{x}^{q}\left( t\right) \right) \left[ 0\right] =\left[
2\right] _{q}x^{k}\left( x-1\right) ^{m}  \label{equation 5}
\end{equation}

By (\ref{equation 100}), we see%
\begin{equation}
\left( D^{N}F_{x}^{q}\left( t\right) \right) \left[ 0\right] =\widetilde{E}%
_{N,q}\left( x\right) \text{ and\ }\left( e^{t}D^{N}F_{x}^{q}\left( t\right)
\right) \left[ 0\right] =\widetilde{E}_{N,q}\left( x\right)
\label{equation 102}
\end{equation}

By expressions of (\ref{equation 5}) and (\ref{equation 102}), we see that%
\begin{equation}
\sum_{j=0}^{\max \left\{ k,m\right\} }\left[ q\binom{k}{j}+\left( -1\right)
^{j}\binom{m}{j}\right] \widetilde{E}_{k+m-j,q}\left( x\right) =\left[ 2%
\right] _{q}x^{k}\left( x-1\right) ^{m}.  \label{equation 6}
\end{equation}

From (\ref{equation 100}), we note that%
\begin{equation}
\frac{d}{dx}\left( \widetilde{E}_{n,q}\left( x\right) \right)
=n\sum_{l=0}^{n-1}\binom{n-1}{l}\widetilde{E}_{l,q}x^{n-1-l}=n\widetilde{E}%
_{n-1,q}\left( x\right)  \label{equation 7}
\end{equation}

By (\ref{equation 7}), we easily see, 
\begin{equation}
\int_{0}^{1}\widetilde{E}_{n,q}\left( x\right) dx=\frac{\widetilde{E}%
_{n+1,q}\left( 1\right) -\widetilde{E}_{n+1,q}}{n+1}=-\frac{\left[ 2\right]
_{q^{-1}}}{n+1}\widetilde{E}_{n+1,q}  \label{equation 8}
\end{equation}

Now, let us consider definition of integral from $0$ to $1$ in (\ref%
{equation 6}), then we have%
\begin{eqnarray}
&&-\left[ 2\right] _{q^{-1}}\sum_{j=0}^{\max \left\{ k,m\right\} }\left[ q%
\binom{k}{j}+\left( -1\right) ^{j}\binom{m}{j}\right] \frac{\widetilde{E}%
_{k+m-j+1,q}}{k+m-j+1}  \label{equation 9} \\
&=&\left[ 2\right] _{q}\left( -1\right) ^{m}B\left( k+1,m+1\right)  \notag \\
&=&\left[ 2\right] _{q}\left( -1\right) ^{m}\frac{\Gamma \left( k+1\right)
\Gamma \left( m+1\right) }{\Gamma \left( k+m+2\right) }  \notag
\end{eqnarray}

where $B\left( m,n\right) $ is beta function which is defined by%
\begin{eqnarray*}
B\left( m,n\right) &=&\int_{0}^{1}x^{m-1}\left( 1-x\right) ^{n-1}dx \\
&=&\frac{\Gamma \left( m\right) \Gamma \left( n\right) }{\Gamma \left(
m+n\right) },\text{ }m>0\text{ and\ }n>0.
\end{eqnarray*}

As a result, we obtain the following theorem

\begin{theorem}
For $n\in 
\mathbb{N}
,$ we have 
\begin{eqnarray*}
&&\sum_{j=1}^{\max \left\{ k,m\right\} }\left[ q\binom{k}{j}+\left(
-1\right) ^{j}\binom{m}{j}\right] \frac{\widetilde{E}_{k+m-j+1,q}}{k+m-j+1}
\\
&=&q\frac{\left( -1\right) ^{m+1}}{\left( k+m+1\right) \binom{k+m}{k}}-\left[
2\right] _{q}\frac{\widetilde{E}_{k+m+1,q}}{k+m+1}.
\end{eqnarray*}
\end{theorem}

Substituting $m=k+1$ into Theorem 1, we readily get 
\begin{eqnarray*}
&&\sum_{j=1}^{k+1}\left[ q\binom{k}{j}+\left( -1\right) ^{j}\binom{k+1}{j}%
\right] \frac{\widetilde{E}_{2k+2-j,q}}{2k+2-j} \\
&=&q\frac{\left( -1\right) ^{k}}{\left( 2k+2\right) \binom{2k+1}{k}}-\left[ 2%
\right] _{q}\frac{\widetilde{E}_{2k+2,q}}{2k+2}.
\end{eqnarray*}

By (\ref{equation 100}), it follows that%
\begin{eqnarray*}
&&\sum_{j=0}^{\max \left\{ k,m\right\} }\left( k+m-j\right) \left[ q\binom{k%
}{j}+\left( -1\right) ^{j}\binom{m}{j}\right] \widetilde{E}%
_{k+m-j-1,q}\left( x\right) \\
&=&\left[ 2\right] _{q}x^{k-1}\left( x-1\right) ^{m-1}\left( \left(
k+m\right) x-k\right) .
\end{eqnarray*}

Let $m=k$ in (\ref{equation 100}), we see that%
\begin{equation*}
\sum_{j=0}^{k}\left[ q\binom{k}{j}+\left( -1\right) ^{j}\binom{k}{j}\right] 
\widetilde{E}_{2k-j,q}\left( x\right) =\left[ 2\right] _{q}x^{k}\left(
x-1\right) ^{k}
\end{equation*}

Last from equality, we discover the following%
\begin{equation}
\left[ 2\right] _{q}\sum_{j=0}^{\left[ \frac{k}{2}\right] }\binom{k}{2j}%
\widetilde{E}_{2k-2j,q}\left( x\right) +\left( q-1\right) \sum_{j=0}^{\left[ 
\frac{k}{2}\right] }\binom{k}{2j+1}\widetilde{E}_{2k-2j-1,q}\left( x\right) =%
\left[ 2\right] _{q}x^{k}\left( x-1\right) ^{k}.  \label{equation 103}
\end{equation}

Here $\left[ .\right] $ is Gauss' symbol. Then, taking integral from $0\ $to 
$1$ both sides of last equality, we get%
\begin{eqnarray*}
&&-\left[ 2\right] _{q^{-1}}\left[ 2\right] _{q}\sum_{j=0}^{\left[ \frac{k}{2%
}\right] }\binom{k}{2j}\frac{\widetilde{E}_{2k-2j+1,q}}{2k-2j+1}+\left[ 2%
\right] _{q^{-1}}\left( 1-q\right) \sum_{j=0}^{\left[ \frac{k}{2}\right] }%
\binom{k}{2j+1}\frac{\widetilde{E}_{2k-2j,q}}{2k-2j} \\
&=&\left[ 2\right] _{q}\left( -1\right) ^{k}B\left( k+1,k+1\right) \\
&=&\frac{\left[ 2\right] _{q}\left( -1\right) ^{k}}{\left( 2k+1\right) 
\binom{2k}{k}}.
\end{eqnarray*}

Consequently, we derive the following theorem

\begin{theorem}
The following identity%
\begin{eqnarray}
&&\left[ 2\right] _{q}\sum_{j=0}^{\left[ \frac{k}{2}\right] }\binom{k}{2j}%
\frac{\widetilde{E}_{2k-2j+1,q}}{2k-2j+1}+\left( q-1\right) \sum_{j=0}^{%
\left[ \frac{k}{2}\right] }\binom{k}{2j+1}\frac{\widetilde{E}_{2k-2j,q}}{%
2k-2j}  \label{equation 10} \\
&=&\frac{q\left( -1\right) ^{k+1}}{\left( 2k+1\right) \binom{2k}{k}}.  \notag
\end{eqnarray}%
is true.
\end{theorem}

In view of (\ref{equation 100}) and (\ref{equation 103}), we discover the
following applications:%
\begin{eqnarray}
&=&\sum_{j=0}^{k+1}\left[ q\binom{k}{j}+\left( -1\right) ^{j}\binom{k+1}{j}%
\right] \widetilde{E}_{2k+1-j,q}\left( x\right)  \label{equation 104} \\
&=&\left[ 2\right] _{q}\widetilde{E}_{2k+1,q}\left( x\right) +\sum_{j=1}^{%
\left[ \frac{k+1}{2}\right] }\left[ q\binom{k}{2j}+\binom{k}{2j}+\binom{k}{%
2j-1}\right] \widetilde{E}_{2k+1-2j,q}\left( x\right)  \notag \\
&&+\sum_{j=0}^{\left[ \frac{k+1}{2}\right] }\left[ q\binom{k}{2j+1}-\binom{k%
}{2j+1}-\binom{k}{2j}\right] \widetilde{E}_{2k-2j,q}\left( x\right)  \notag
\\
&=&-\left[ \sum_{j=0}^{\left[ \frac{k}{2}\right] }\binom{k}{2j}\widetilde{E}%
_{2k-2j,q}\left( x\right) +\frac{q-1}{1+q}\sum_{j=0}^{\left[ \frac{k}{2}%
\right] }\binom{k}{2j+1}\widetilde{E}_{2k-2j+1}\left( x\right) \right] 
\notag \\
&&+\left[ 2\right] _{q}\sum_{j=0}^{\left[ \frac{k}{2}\right] }\binom{k}{2j}%
\widetilde{E}_{2k+1-2j,q}\left( x\right) +\sum_{j=1}^{\left[ \frac{k}{2}%
\right] }\binom{k}{2j-1}\widetilde{E}_{2k+1-2j,q}\left( x\right)  \notag \\
&&+\left( q-1\right) \sum_{j=0}^{\left[ \frac{k}{2}\right] }\binom{k}{2j+1}%
\widetilde{E}_{2k-2j,q}\left( x\right) +\frac{q-1}{1+q}\sum_{j=0}^{\left[ 
\frac{k}{2}\right] }\binom{k}{2j+1}\widetilde{E}_{2k-2j+1}\left( x\right) 
\notag
\end{eqnarray}

By expressions (\ref{equation 103}) and (\ref{equation 104}), we have the
following Theorem

\begin{theorem}
For $k\in 
\mathbb{N}
$, we have%
\begin{eqnarray}
&&\left[ 2\right] _{q}\sum_{j=0}^{\left[ \frac{k}{2}\right] }\binom{k}{2j}%
\widetilde{E}_{2k+1-2j,q}\left( x\right) +\sum_{j=1}^{\left[ \frac{k}{2}%
\right] }\binom{k}{2j-1}\widetilde{E}_{2k+1-2j,q}\left( x\right)
\label{equation 105} \\
&&+\left( q-1\right) \sum_{j=0}^{\left[ \frac{k}{2}\right] }\binom{k}{2j+1}%
\left[ \widetilde{E}_{2k-2j,q}\left( x\right) +\frac{1}{1+q}\widetilde{E}%
_{2k-2j+1}\left( x\right) \right]  \notag \\
&=&x^{k}\left( x-1\right) ^{k}\left( \left[ 2\right] _{q}x-q\right)  \notag
\end{eqnarray}
\end{theorem}

\section{$p$\textbf{-adic integral on\ }$%
\mathbb{Z}
_{p}$\textbf{\ associated with Kim's }$q$\textbf{-Euler polynomials}}

In this section, we consider Kim's $q$-Euler polynomials by means of $p$%
-adic $q$-integral on $%
\mathbb{Z}
_{p}$. Now we start with the following assertion.

Let $m,k\in 
\mathbb{N}
$, Then by (\ref{equation 6}), 
\begin{eqnarray*}
I_{1} &=&\left[ 2\right] _{q}\int_{%
\mathbb{Z}
_{p}}x^{k}\left( x-1\right) ^{m}d\mu _{-q}\left( x\right) \\
&=&\left[ 2\right] _{q}\sum_{l=0}^{m}\binom{m}{l}\left( -1\right)
^{m-l}\int_{%
\mathbb{Z}
_{p}}x^{l+k}d\mu _{-q}\left( x\right) \\
&=&\left[ 2\right] _{q}\sum_{l=0}^{m}\binom{m}{l}\left( -1\right) ^{m-l}%
\widetilde{E}_{l+k,q}
\end{eqnarray*}

On the other hand, right hand side of (\ref{equation 6}), 
\begin{eqnarray*}
I_{2} &&=\sum_{j=0}^{\max \left\{ k,m\right\} }\left[ q\binom{k}{j}+\left(
-1\right) ^{j}\binom{m}{j}\right] \sum_{l=0}^{k+m-j}\binom{k+m-j}{l}%
\widetilde{E}_{k+m-j-l,q}\int_{%
\mathbb{Z}
_{p}}x^{l}d\mu _{-q}\left( x\right) \\
&=&\sum_{j=0}^{\max \left\{ k,m\right\} }\left[ q\binom{k}{j}+\left(
-1\right) ^{j}\binom{m}{j}\right] \sum_{l=0}^{k+m-j}\binom{k+m-j}{l}%
\widetilde{E}_{k+m-j-l,q}\widetilde{E}_{l,q}
\end{eqnarray*}

Equating $I_{1}$ and $I_{2}$, we get the following theorem

\begin{theorem}
For $m,k\in 
\mathbb{N}
$, we have%
\begin{eqnarray*}
&&\sum_{j=0}^{\max \left\{ k,m\right\} }\left[ q\binom{k}{j}+\left(
-1\right) ^{j}\binom{m}{j}\right] \sum_{l=0}^{k+m-j}\binom{k+m-j}{l}%
\widetilde{E}_{k+m-j-l,q}\widetilde{E}_{l,q} \\
&=&\left[ 2\right] _{q}\sum_{l=0}^{m}\binom{m}{l}\left( -1\right) ^{m-l}%
\widetilde{E}_{l+k,q}.
\end{eqnarray*}
\end{theorem}

Let us take fermionic $p$-adic $q$-inetgral on $%
\mathbb{Z}
_{p}$ left hand side of (\ref{equation 105}), we get%
\begin{eqnarray*}
I_{3} &=&\int_{%
\mathbb{Z}
_{p}}x^{k}\left( x-1\right) ^{k}\left( \left[ 2\right] _{q}x-q\right) d\mu
_{-q}\left( x\right) \\
&=&\left[ 2\right] _{q}\sum_{l=0}^{k}\binom{k}{l}\left( -1\right)
^{k-l}\int_{%
\mathbb{Z}
_{p}}x^{k+l+1}d\mu _{-q}\left( x\right) -q\sum_{l=0}^{k}\binom{k}{l}\left(
-1\right) ^{k-l}\int_{%
\mathbb{Z}
_{p}}x^{k+l}d\mu _{-q}\left( x\right) \\
&=&\left[ 2\right] _{q}\sum_{l=0}^{k}\binom{k}{l}\left( -1\right) ^{k-l}%
\widetilde{E}_{k+l+1,q}-q\sum_{l=0}^{k}\binom{k}{l}\left( -1\right) ^{k-l}%
\widetilde{E}_{k+l,q}
\end{eqnarray*}

In other word, we consider right hand side of (\ref{equation 105}) as
follows:%
\begin{eqnarray*}
I_{4} &=&\left[ 2\right] _{q}\sum_{j=0}^{\left[ \frac{k}{2}\right] }\binom{k%
}{2j}\sum_{l=0}^{2k-2j+1}\binom{2k-2j+1}{l}\widetilde{E}_{2k+1-2j-l,q}\int_{%
\mathbb{Z}
_{p}}x^{l}d\mu _{-q}\left( x\right) \\
&&+\sum_{j=1}^{\left[ \frac{k}{2}\right] }\binom{k}{2j-1}\sum_{l=0}^{2k-2j+1}%
\binom{2k-2j+1}{l}\widetilde{E}_{2k+1-2j-l,q}\int_{%
\mathbb{Z}
_{p}}x^{l}d\mu _{-q}\left( x\right) \\
&&+\sum_{j=0}^{\left[ \frac{k}{2}\right] }\binom{k}{2j+1}\left[ 
\begin{array}{c}
\left( q-1\right) \sum_{j=0}^{2k-2j}\binom{2k-2j}{l}\widetilde{E}%
_{2k-2j-l,q}\int_{%
\mathbb{Z}
_{p}}x^{l}d\mu _{-q}\left( x\right) \\ 
+\frac{q-1}{1+q}\sum_{l=0}^{2k-2j+1}\binom{2k-2j+1}{l}\widetilde{E}%
_{2k-2j-l+1}\int_{%
\mathbb{Z}
_{p}}x^{l}d\mu _{-q}\left( x\right)%
\end{array}%
\right] \\
&=&\left[ 2\right] _{q}\sum_{j=0}^{\left[ \frac{k}{2}\right] }\binom{k}{2j}%
\sum_{l=0}^{2k-2j+1}\binom{2k-2j+1}{l}\widetilde{E}_{2k+1-2j-l,q}\widetilde{E%
}_{l,q} \\
&&+\sum_{j=1}^{\left[ \frac{k}{2}\right] }\binom{k}{2j-1}\sum_{l=0}^{2k-2j+1}%
\binom{2k-2j+1}{l}\widetilde{E}_{2k+1-2j-l,q}\widetilde{E}_{l,q} \\
&&+\sum_{j=0}^{\left[ \frac{k}{2}\right] }\binom{k}{2j+1}\left[ 
\begin{array}{c}
\left( q-1\right) \sum_{j=0}^{2k-2j}\binom{2k-2j}{l}\widetilde{E}_{2k-2j-l,q}%
\widetilde{E}_{l,q} \\ 
+\frac{q-1}{1+q}\sum_{l=0}^{2k-2j+1}\binom{2k-2j+1}{l}\widetilde{E}%
_{2k-2j-l+1}\widetilde{E}_{l,q}%
\end{array}%
\right]
\end{eqnarray*}

Equating $I_{3}$ and $I_{4}$, we get the following theorem

\begin{theorem}
For $k\in 
\mathbb{N}
$, we have%
\begin{eqnarray*}
&&\sum_{l=0}^{k}\binom{k}{l}\left( -1\right) ^{k-l}\left[ \left[ 2\right]
_{q}\widetilde{E}_{k+l+1,q}-q\widetilde{E}_{k+l,q}\right] \\
&=&\left[ 2\right] _{q}\sum_{j=0}^{\left[ \frac{k}{2}\right] }\binom{k}{2j}%
\sum_{l=0}^{2k-2j+1}\binom{2k-2j+1}{l}\widetilde{E}_{2k+1-2j-l,q}\widetilde{E%
}_{l,q} \\
&&+\sum_{j=1}^{\left[ \frac{k}{2}\right] }\binom{k}{2j-1}\sum_{l=0}^{2k-2j+1}%
\binom{2k-2j+1}{l}\widetilde{E}_{2k+1-2j-l,q}\widetilde{E}_{l,q} \\
&&+\sum_{j=0}^{\left[ \frac{k}{2}\right] }\binom{k}{2j+1}\left\{ 
\begin{array}{c}
\left( q-1\right) \sum_{j=0}^{2k-2j}\binom{2k-2j}{l}\widetilde{E}_{2k-2j-l,q}%
\widetilde{E}_{l,q} \\ 
+\frac{q-1}{1+q}\sum_{l=0}^{2k-2j+1}\binom{2k-2j+1}{l}\widetilde{E}%
_{2k-2j-l+1}\widetilde{E}_{l,q}%
\end{array}%
\right\}
\end{eqnarray*}
\end{theorem}

Now, we consider (\ref{equation 6}) and (\ref{equation 100}) by means of $q$%
-Volkenborn integral. Then, by (\ref{equation 6}), we see%
\begin{eqnarray*}
&&\left[ 2\right] _{q}\int_{%
\mathbb{Z}
_{p}}x^{k}\left( x-1\right) ^{m}d\mu _{q}\left( x\right) \\
&=&\left[ 2\right] _{q}\sum_{l=0}^{m}\binom{m}{l}\left( -1\right)
^{m-l}\int_{%
\mathbb{Z}
_{p}}x^{l+k}d\mu _{q}\left( x\right) \\
&=&\left[ 2\right] _{q}\sum_{l=0}^{m}\binom{m}{l}\left( -1\right) ^{m-l}%
\widetilde{B}_{l+k,q}
\end{eqnarray*}

On the other hand,%
\begin{eqnarray*}
&&\sum_{j=0}^{\max \left\{ k,m\right\} }\left[ q\binom{k}{j}+\left(
-1\right) ^{j}\binom{m}{j}\right] \sum_{l=0}^{k+m-j}\binom{k+m-j}{l}%
\widetilde{E}_{k+m-j-l,q}\int_{%
\mathbb{Z}
_{p}}x^{l}d\mu _{q}\left( x\right) \\
&=&\sum_{j=0}^{\max \left\{ k,m\right\} }\left[ q\binom{k}{j}+\left(
-1\right) ^{j}\binom{m}{j}\right] \sum_{l=0}^{k+m-j}\binom{k+m-j}{l}%
\widetilde{E}_{k+m-j-l,q}\widetilde{B}_{l,q}
\end{eqnarray*}

Therefore, we get the following theorem

\begin{theorem}
For $m,k\in 
\mathbb{N}
$, we have%
\begin{eqnarray*}
&&\left[ 2\right] _{q}\sum_{l=0}^{m}\binom{m}{l}\left( -1\right) ^{m-l}%
\widetilde{B}_{l+k,q} \\
&=&\sum_{j=0}^{\max \left\{ k,m\right\} }\left[ q\binom{k}{j}+\left(
-1\right) ^{j}\binom{m}{j}\right] \sum_{l=0}^{k+m-j}\binom{k+m-j}{l}%
\widetilde{E}_{k+m-j-l,q}\widetilde{B}_{l,q}
\end{eqnarray*}
\end{theorem}

By using fermionic $p$-adic $q$-integral on $%
\mathbb{Z}
_{p}$ left hand side of (\ref{equation 105}), we get%
\begin{eqnarray*}
I_{5} &=&\left[ 2\right] _{q}\int_{%
\mathbb{Z}
_{p}}x^{k}\left( x-1\right) ^{k}\left( \left[ 2\right] x-q\right) d\mu
_{q}\left( x\right) \\
&=&\left[ 2\right] _{q}\sum_{l=0}^{k}\binom{k}{l}\left( -1\right)
^{k-l}\int_{%
\mathbb{Z}
_{p}}x^{k+l+1}d\mu _{q}\left( x\right) -q\sum_{l=0}^{k}\binom{k}{l}\left(
-1\right) ^{k-l}\int_{%
\mathbb{Z}
_{p}}x^{k+l}d\mu _{q}\left( x\right) \\
&=&\left[ 2\right] _{q}\sum_{l=0}^{k}\binom{k}{l}\left( -1\right) ^{k-l}%
\widetilde{B}_{k+l+1,q}-q\sum_{l=0}^{k}\binom{k}{l}\left( -1\right) ^{k-l}%
\widetilde{B}_{k+l,q}
\end{eqnarray*}

Also, we consider right hand side of (\ref{equation 105}) as follows:%
\begin{eqnarray*}
I_{6} &=&\left[ 2\right] _{q}\sum_{j=0}^{\left[ \frac{k}{2}\right] }\binom{k%
}{2j}\sum_{l=0}^{2k-2j+1}\binom{2k-2j+1}{l}\widetilde{E}_{2k+1-2j-l,q}\int_{%
\mathbb{Z}
_{p}}x^{l}d\mu _{q}\left( x\right) \\
&&+\sum_{j=1}^{\left[ \frac{k}{2}\right] }\binom{k}{2j-1}\sum_{l=0}^{2k-2j+1}%
\binom{2k-2j+1}{l}\widetilde{E}_{2k+1-2j-l,q}\int_{%
\mathbb{Z}
_{p}}x^{l}d\mu _{q}\left( x\right) \\
&&+\sum_{j=0}^{\left[ \frac{k}{2}\right] }\binom{k}{2j+1}\left[ 
\begin{array}{c}
\left( q-1\right) \sum_{j=0}^{2k-2j}\binom{2k-2j}{l}\widetilde{E}%
_{2k-2j-l,q}\int_{%
\mathbb{Z}
_{p}}x^{l}d\mu _{q}\left( x\right) \\ 
+\frac{q-1}{1+q}\sum_{l=0}^{2k-2j+1}\binom{2k-2j+1}{l}\widetilde{E}%
_{2k-2j-l+1}\int_{%
\mathbb{Z}
_{p}}x^{l}d\mu _{q}\left( x\right)%
\end{array}%
\right] \\
&=&\left[ 2\right] _{q}\sum_{j=0}^{\left[ \frac{k}{2}\right] }\binom{k}{2j}%
\sum_{l=0}^{2k-2j+1}\binom{2k-2j+1}{l}\widetilde{E}_{2k+1-2j-l,q}\widetilde{B%
}_{l,q} \\
&&+\sum_{j=1}^{\left[ \frac{k}{2}\right] }\binom{k}{2j-1}\sum_{l=0}^{2k-2j+1}%
\binom{2k-2j+1}{l}\widetilde{E}_{2k+1-2j-l,q}\widetilde{B}_{l,q} \\
&&+\sum_{j=0}^{\left[ \frac{k}{2}\right] }\binom{k}{2j+1}\left[ 
\begin{array}{c}
\left( q-1\right) \sum_{j=0}^{2k-2j}\binom{2k-2j}{l}\widetilde{E}_{2k-2j-l,q}%
\widetilde{B}_{l,q} \\ 
+\frac{q-1}{1+q}\sum_{l=0}^{2k-2j+1}\binom{2k-2j+1}{l}\widetilde{E}%
_{2k-2j-l+1}\widetilde{B}_{l,q}%
\end{array}%
\right]
\end{eqnarray*}

Equating $I_{5}$ and $I_{6}$, we get the following Corollary

\begin{corollary}
For $k\in 
\mathbb{N}
$, we get%
\begin{eqnarray*}
&&\sum_{l=0}^{k}\binom{k}{l}\left( -1\right) ^{k-l}\left[ \left[ 2\right]
_{q}\widetilde{B}_{k+l+1,q}-q\widetilde{B}_{k+l,q}\right] \\
&=&\left[ 2\right] _{q}\sum_{j=0}^{\left[ \frac{k}{2}\right] }\binom{k}{2j}%
\sum_{l=0}^{2k-2j+1}\binom{2k-2j+1}{l}\widetilde{E}_{2k+1-2j-l,q}\widetilde{B%
}_{l,q} \\
&&+\sum_{j=1}^{\left[ \frac{k}{2}\right] }\binom{k}{2j-1}\sum_{l=0}^{2k-2j+1}%
\binom{2k-2j+1}{l}\widetilde{E}_{2k+1-2j-l,q}\widetilde{B}_{l,q} \\
&&+\sum_{j=0}^{\left[ \frac{k}{2}\right] }\binom{k}{2j+1}\left\{ 
\begin{array}{c}
\left( q-1\right) \sum_{j=0}^{2k-2j}\binom{2k-2j}{l}\widetilde{E}_{2k-2j-l,q}%
\widetilde{B}_{l,q} \\ 
+\frac{q-1}{1+q}\sum_{l=0}^{2k-2j+1}\binom{2k-2j+1}{l}\widetilde{E}%
_{2k-2j-l+1}\widetilde{B}_{l,q}%
\end{array}%
\right\}
\end{eqnarray*}%
\ 
\end{corollary}

\end{document}